%% file: ms.tex
\documentclass[letterpaper, 10 pt, conference]{ieeeconf}

\IEEEoverridecommandlockouts                              

\input{_paper}
\usepackage[compress, sort]{cite}

\title{\LARGE \bf
Hierarchical distributed scenario-based model predictive control \mbox{of interconnected microgrids}*
}

\author{T. Alissa Schenck$^{1}$ and Christian A. Hans$^{2}$
\thanks{*Note that this work is based on the master thesis \enquote{Stochastic model predictive control of interconnected microgrids} by T. Alissa Schenck.}
\thanks{$^{1}$T. Alissa Schenck received a master's degree from
	Technische Universität Berlin,
	Germany,
        {\tt \href{mailto:ta-schenck@posteo.de}{ta-schenck@posteo.de}}.}%
\thanks{$^{2}$Christian A. Hans is with the
	Automation and Sensorics in Networked Systems Group,
	University of Kassel,
	Germany,
	{\tt \href{mailto:hans@uni-kassel.de}{hans@uni-kassel.de}}.}%
}

\begin{document}

\maketitle

\begin{abstract}
Microgrids are autonomous clusters of generators, storage units and loads.
Special requirements arise in interconnected operation:
control schemes that do not require individual microgrids to disclose information about their internal structure and operating objectives are preferred for privacy reasons.
Moreover, a safe and economically meaningful operation shall be achieved in presence of uncertain load and weather-dependent availability of renewable infeed.
In this paper, we propose a hierarchical distributed model predictive control approach that satisfies these requirements.
Specifically, we demonstrate that costs and safety of supply can be improved through a scenario-based stochastic control scheme.
In a numerical case study, our approach is compared to a certainty equivalence and a prescient scheme.
The results illustrate good performance as well as sufficiently fast convergence.
\end{abstract}

\input{introduction}
\input{gridModel}
\input{problemFormulation}
\input{distributedSolution}

\input{caseStudy}
\input{conclusions}

\bibliographystyle{IEEEtran}
\bibliography{literature}

\end{document}

%% file: _paper.tex
\usepackage[utf8]{inputenc}
\usepackage[english]{babel}

\IEEEoverridecommandlockouts 

\usepackage{csquotes}

\usepackage{amsmath,amssymb,amsfonts}
\usepackage{algorithm}
\usepackage{graphicx}
\usepackage{textcomp}
\usepackage{xcolor}

\definecolor{vgRed}{RGB}{193, 48, 24}
\definecolor{vgOrange}{RGB}{243, 111, 19}
\definecolor{vgYellow}{RGB}{235, 203, 56}
\definecolor{vgGreen}{RGB}{162, 185, 105}
\definecolor{vgLightBlue}{RGB}{13, 149, 188}
\definecolor{vgDarkBlue}{RGB}{6, 56, 81}

\colorlet{mycolor1}{black}
\colorlet{mycolor2}{black}
\colorlet{mycolor3}{black}

\usepackage{orcidlink}
\usepackage[capitalise]{cleveref}
\crefname{equation}{}{}
\Crefname{equation}{Equation}{Equations}

\crefname{figure}{Figure}{Figures}
\crefname{section}{Section}{Sections}
\crefname{algorithm}{Algorithm}{Algorithms}
\crefname{problem}{Problem}{Problems}
\crefname{proposition}{Proposition}{Propositions}
\crefname{remark}{Remark}{Remarks}
\Crefname{remark}{Remark}{Remarks}

\newtheorem{remark}{Remark}
\newtheorem{problem}{Problem}

\usepackage{gitinfo2}

\usepackage{acronym}
\acrodef{admm}[ADMM]{alternating direction method of multipliers}
\acrodef{arima}[ARIMA]{autoregressive integrated moving average}
\acrodef{der}[DER]{distributed energy ressources}
\acrodef{dnp3}[DNP3]{distributed network protocol}
\acrodef{epa}[EPA]{enhanced performance architecture}
\acrodef{ghg}[GHG]{greenhouse gas}
\acrodef{lan}[LAN]{local area network}
\acrodef{mg}[MG]{microgrid}
\acrodef{mpc}[MPC]{model predictive control}
\acrodef{pcc}[PCC]{point of common coupling}
\acrodef{res}[RES]{renewable energy source}
\acrodef{sarima}[SARIMA]{seasonal autoregressive integrated moving average}
\acrodef{tcp/ip}[TCP/IP]{transmission control protocol / internet protocol}
\acrodef{wan}[WAN]{wide area network}
\acrodef{ac}[AC]{alternating current} \acused{ac}
\acrodef{dc}[DC]{direct current} \acused{dc}
\acrodef{mip}[MIP]{mixed-integer problem} 
\acrodef{miqp}[MIQP]{mixed-integer quadratic problem}

\usepackage{booktabs}
\usepackage{multirow}
\usepackage{units}
\let\labelindent\relax
\usepackage[inline, shortlabels]{enumitem}
\usepackage{csvsimple}

\usepackage{pgfplots}
\usetikzlibrary{positioning}
\usetikzlibrary{shapes.multipart}
\usetikzlibrary{arrows}
\usetikzlibrary{decorations.pathreplacing,decorations.markings}
\usetikzlibrary{calc}

\pgfplotsset{compat=1.18}
\usepgfplotslibrary{external}
\tikzset{external/system call={pdflatex \tikzexternalcheckshellescape -halt-on-error
		-interaction=batchmode -jobname "\image" "\texsource"}}
\newcommand{\includetikz}[1]{%
    \includegraphics{tikzPictures/#1.pdf}%
}

\usepackage{todonotes}
\makeatletter
\renewcommand{\todo}[2][]{\tikzexternaldisable\@todo[#1]{#2}\tikzexternalenable}
\renewcommand{\missingfigure}[2][]{\tikzexternaldisable\@missingfigure[#1]{#2}\tikzexternalenable}
\newcommand{\atAj}[1]{\tikzexternaldisable\@todo[inline, color=blue!20]{@Ajay: #1}\tikzexternalenable}
\newcommand{\atCh}[1]{\tikzexternaldisable\@todo[inline, color=red!20]{@Chris: #1}\tikzexternalenable}
\newcommand{\fromAj}[1]{\tikzexternaldisable\@todo[inline, color=red!20]{Ajay says: #1}\tikzexternalenable}
\makeatother

\def\R{\mathbb{R}}

\def\N{\mathbb{N}}

\newcommand{\tran}{^{\mathstrut\scriptscriptstyle\top}}
\newcommand{\itran}{{\mathstrut\scriptscriptstyle\top}}

\newcommand{\diag}{\operatorname{diag}}

\newcommand{\ivar}[3]{#1_{#2}^{#3}}
\newcommand{\norm}[1]{\| #1 \|}
\newcommand{\stage}[2][i]{\operatorname{stage}_{#1}(#2)}
\newcommand{\nodes}[2][i]{\operatorname{nodes}_{#1}(#2)}
\newcommand{\child}[2][i]{\operatorname{child}_{#1}(#2)}
\newcommand{\anc}[2][i]{\operatorname{anc}_{#1}(#2)}
\newcommand{\smallSum}{\textstyle\sum}

\DeclareMathOperator*{\minimize}{\operatorname{minimize}}

\DeclareMathOperator*{\argmin}{\ensuremath{\operatorname{arg\,min}}}
\newcommand{\qed}{\hfill $\bullet$}

\makeatletter
\newcommand{\subalign}[1]{%
  \vcenter{%
    \Let@ \restore@math@cr \default@tag
    \baselineskip\fontdimen10 \scriptfont\tw@
    \advance\baselineskip\fontdimen12 \scriptfont\tw@
    \lineskip\thr@@\fontdimen8 \scriptfont\thr@@
    \lineskiplimit\lineskip
    \ialign{\hfil$\m@th\scriptstyle##$&$\m@th\scriptstyle{}##$\hfil\crcr
      #1\crcr
    }%
  }%
}
\makeatother

\pgfplotsset{every axis/.append style={semithick,tick style={major tick
            length=4pt,semithick,black}}}

\pgfkeys{/pgfplots/x axis shift down/.style={
        x axis line style={yshift=-#1},
        xtick style={yshift=-#1},
        tick align=outside,
        xticklabel shift={#1}}}

\pgfkeys{/pgfplots/y axis shift left/.style={
        y axis line style={xshift=-#1},
        ytick style={xshift=-#1},
        yticklabel shift={#1},
        scaled y ticks = false,
        tick align=outside,
        y tick label style={/pgf/number format/fixed,
        }
    }
}

\pgfplotsset{myPlot/.style={%
        line width = 0.7pt,
        separate axis lines,
        axis x line*=bottom,
        x axis shift down = 3pt,
        enlarge x limits=false,
        axis y line*=left,
        y axis shift left = 6pt,
        enlarge y limits={abs=.25pt},
        enlarge x limits={abs=.25pt},
    }
}

%% file: introduction.tex

\section{Introduction}
\label{sec:introduction}

Climate change is one of the major challenges of the 21st century.
\Acp{res} play a central role in tackling this challenge.
As infeed of \acp{res} increases, new questions arise:
how to ensure a reliable operation in presence of uncertain \acp{res}
and
how to keep the complexity manageable despite a large number small-scale \acp{res}?

\Acp{mg} have emerged as a promising approach to answer these questions.
An \ac{mg} clusters a collection of loads and distributed units \cite{OME+2014}, such as storage units, \acp{res} and conventional generators.
Internally, it is controlled in a way that makes it appear as a single subsystem \cite{PL2006} to the outside world which acts as a load or a generator \cite{OME+2014} when connected to a utility grid.
Compared to islanded operation, trading energy between interconnected \acp{mg} allows to further increase infeed of \acp{res} with different generation patterns.

Operation control, also referred to as energy management, aims to provide power setpoints to the units and thereby control the energy of storages.
For this task, \ac{mpc} has been widely employed.
Here, optimization problems are solved to find control actions that minimize an objective subject to constrains which model the system behavior and account for limits.
Decisions have to be made in presence of uncertain load and renewable infeed.
In state-of-the-art approaches, nominal forecasts are often assumed to be certain which holds the risk of constraint violations in presence of prediction errors.
Stochastic approaches, that employ forecast probability distributions, on the other hand, allow to increase robustness and performance \cite{Han2021}.

Stochastic \ac{mpc} of \acp{mg} was actively explored in recent years.
In \cite{MSM2014}, a scenario-based method for the optimal operation of a single grid-connected \ac{mg} was proposed.
Kou et al. consider a similar setup in \cite{KLG2018}.
In \cite{PVC+2016}, a scenario-based approach consisting of optimal generation scheduling and \ac{mpc} is presented.
The authors of \cite{MD2017} propose a scenario-based approach that employs heuristics for optimal control.
Heymann et al. model the demand dynamics with stochastic differential equation and solve a management problem \cite{HBSJ2016}.
Only few have studied stochastic control of networks of interconnected \acp{mg}, probably because such problems can be hard to decompose and stochastic \acp{mip} often do not scale well.
To our knowledge, only in \cite{BTAG2019}, a multi-\ac{mg} system with uncertain \acp{res} and load as well as power exchange between \acp{mg} is considered.
However, the approach relies on day-ahead scheduling which can be inflexible in grid with high share of uncertain \aclp{res}.


In this work, we employ a hierarchical distributed \ac{mpc} approach \cite{HBR+2018} in a setting where scenario-based stochastic local problems along the lines of \cite{HSB+2015} are considered.
In this context, the following contributions are made.
\begin{enumerate*}[(i)]
	\item
		Uncertainties in renewable generation and loads are handled by inclusion of forecast probability distributions.
		In the overall \ac{mpc} scheme, we assume that fluctuations are covered within each \ac{mg} and power flow over the \ac{ac} grid is certain.
	\item
		We solve the \ac{mpc} problem through a distributed algorithm that is based on the widely employed \ac{admm} (see, e.g., \cite{HBR+2018,LFSS2020,RZM+2019}).
		Control actions are found by alternately solving
			local stochastic optimization problems at individual \acp{mg}
			and
			a certainty equivalence optimization problem at a central coordinator that takes care of the power flow between the \acp{mg}.
		The algorithm allows to respect the privacy and independence of individual \acp{mg}.
		In addition, the computational complexity scales well with the number of \acp{mg} such that solve times remain manageable.
	\item
		In a comprehensive case study, the stochastic approach is thoroughly compared with the certainty equivalence approach from \cite{HBR+2018} and a prescient \ac{mpc}.
		Several closed-loop performance metrics as well as solve times and number of iterations are assessed.
		The results highlight the increased security of the novel approach compared to certainty equivalence \ac{mpc}.
\end{enumerate*}


The remainder of this paper is structured as follows.
In \cref{sec:gridModel}, the plant model and the model of the uncertain forecast are introduced.
In \cref{sec:problemFormulation}, a problem formulation is provided and in \cref{sec:distributedSolution}, a distributed solution proposed.
Finally, in \cref{sec:caseStudy} simulation results are discussed.

\subsection{Notation}
\label{sec:notation}

The set of positive integers is $\mathbb{N}$ and the set of nonnegative integers $\mathbb{N}_0$.
The set of Boolean variables is $\mathbb{B} = \{ 1, 0 \}$.
The set $\{ x \vert x \in \mathbb{N}_0 \wedge a \leq x \leq b \}$ is shortly referred to as $\mathbb{N}_{[a,b]}$.
The set of real numbers is $\mathbb{R}$,
the set of negative real numbers $\mathbb{R}_{< 0}$ and
the set of positive real numbers $\mathbb{R}_{> 0}$.
Likewise, $\mathbb{R}_{\leq 0}$ denotes nonpositive and $\mathbb{R}_{\geq 0}$ nonnegative real numbers.

Let $\mathbf{1}_n$ be the $n$-dimensional column vector of all ones.
Moreover, $\mathbf{0}_{m \times n}$ is the ${m \times n}$ matrix of all zeros and
$\mathbf{I}_n$ the $n \times n$ identity matrix.
Consider a vector $a = [a_1 ~ \cdots ~ a_n]\tran$.
Then $\diag(a)$ is a diagonal matrix with entries $a_i$, ${i\in\mathbb{N}_{[1,n]}}$.
Finally, $\norm{a}_2$ is the Euclidean norm of vector $a$.

%% file: gridModel.tex

\section{Model}
\label{sec:gridModel}

\begin{figure}
	\vspace{0.5em}
	\centering
	\tikzset{external/export next=true}
	\includetikz{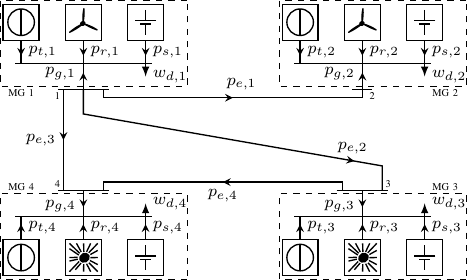}
	\caption{Topology of four interconnected \ac{ac} \acp{mg} from \cite{HBR+2018}.}
	\label{fig:interconnectedMgScheme}
\end{figure}

In what follows, the model of a network of \ac{ac} \acp{mg} (see, e.g., \cref{fig:interconnectedMgScheme}) is derived.
Moreover, the representation of uncertain forecasts in the form of scenario trees is discussed.

\subsection{\acs{mg} model variables}
\label{sec:mg-notation}

\begin{table}
	\vspace{0.7em}
	\centering
	\caption{Variables of \ac{mg} $i$ by unit type}
	\label{tab:domains}
	\begin{tabular}{lllc}
		\toprule
			 & Variable & Type & Symbol \\
		\midrule
		Convent. & Setpoint & Control input & $u_{t,i} \in \mathbb{R}_{\geq 0}^{T_i}$ \\[0.3em]
				 & On/off switch & Control input & $\delta_{t,i} \in \mathbb{B}^{T_i}$ \\[0.3em]
				 & Power & Auxiliary & $p_{t,i} \in \mathbb{R}_{\geq 0}^{T_i}$ \\
		\midrule
		Storage  & Setpoint & Control input & $u_{s,i} \in \mathbb{R}^{S_i}$ \\[0.3em]
				 & Energy & State & $x_i \in \mathbb{R}_{\geq 0}^{S_i}$ \\[0.3em]
				 & Power & Auxiliary & $p_{s,i} \in \mathbb{R}^{S_i}$ \\[0.3em]
				 & Slack var. & Auxiliary & $\sigma_i \in \mathbb{R}_{\geq 0}^{S_i}$ \\
		\midrule
		\ac{res} & Setpoint & Control input & $u_{r,i} \in \mathbb{R}_{\geq 0}^{R_i}$ \\[0.3em]
				 & Weather-dependent & Uncertain input & $w_{r,i} \in \mathbb{R}_{\geq 0}^{R_i}$\\[0.0em]
				 & available power \\
				 & Power & Auxiliary & $p_{r,i} \in \mathbb{R}_{\geq 0}^{R_i}$ \\[0.3em]
				 & Boolean var. for & Auxiliary & $\delta_{r,i} \in \mathbb{B}^{R_i}$\\[0.0em]
				 & $\min$-operator \\[0.0em]
		\midrule
		Load 	 & Power & Uncertain input & $w_{d,i} \in \mathbb{R}_{\leq 0}^{D_i}$\\
			\midrule
		Grid & \acs{pcc} Power & Auxiliary & $p_{g,i} \in \mathbb{R}$ \\[0.3em]
			 & Power sharing var. & Auxiliary & $\rho_{i} \in \mathbb{R}$ \\
		\bottomrule
	\end{tabular}
\end{table}

Consider $I\in\N$ \acp{mg} which are uniquely indexed by elements of the set $\mathbb{I} = \mathbb{N}_{[1, I]}$.
Each \ac{mg} $i\in\mathbb{I}$ is composed of
$T_i \in \mathbb{N}$ conventional units,
$S_i \in \mathbb{N}$ storage units,
$R_i \in \mathbb{N}$ \acp{res} and
$D_i \in \mathbb{N}$ loads.
The variables associated with \ac{mg}~$i$ are collected in \cref{tab:domains}.
For the units, active sign convention is used, i.e., positive values indicate that power is provided, negative values that it is consumed.
For \ac{mg}~$i$, the power setpoints of all units are collected in $u_i = [ \ivar{u}{t,i}{\itran} ~ \ivar{u}{s,i}{\itran} ~ \ivar{u}{r,i}{\itran}]\tran$ and
the control inputs combined in $v_i=[ u_i\tran ~ \delta_{t,i}\tran ]\tran$.
The power values are collected in $p_i = [ \ivar{p}{t,i}{\itran} ~ \ivar{p}{s,i}{\itran} ~ \ivar{p}{r,i}{\itran}]\tran$ which is used to form $q_i=[ p_i\tran ~ \delta_{r,i}\tran ~ \sigma_{i}\tran ~ \rho_i ]\tran$.
The state $x_i$ represents the energy that is contained in the storage units.
Finally, $w_i=[ w_{r,i}\tran ~ w_{d,i}\tran ]\tran$ is the vector of uncertain inputs.

\subsection{Uncertainty model}
\label{sec:uncertaintyModel}

The forecasts of uncertain load and weather-dependent available renewable infeed are modeled using scenario trees.
This section provides a brief introduction on their generation and structure.
Note that it is strongly motivated by \cite{Han2021}.

\subsubsection{From forecast scenarios to scenario trees}

Consider a finite number of independent equiprobable forecast scenarios proceeding into the future for $J$ prediction steps.
For a large number of scenarios, the underlying probability distribution is accurately approximated.
Since a large number of scenarios leads to complex optimization problems, more compact representations of probability distributions, such as scenario trees, are desirable.
In this paper, forecast scenarios of available wind power are derived using \ac{arima} models.
The available photovoltaic power and load demand forecast scenarios are generated using seasonal \ac{arima} models.
On these scenarios, forward selection \cite{HR2003} is applied to obtain scenario-trees.

\subsubsection{Scenario tree structure}
\label{sec:scenario-structure}

Scenario trees are collections of ${M_i \in \mathbb{N}}$ nodes which are partitioned into $J \in \mathbb{N}$ stages, i.e., prediction steps $j \in \mathbb{N}_{[0,J]}$.
Let us collect the unique indices of all nodes in $\mathbb{M}_i = \mathbb{N}_{[0, M_i-1]}$.
The $\stage{m}$ operator provides the stage of node ${m\in\mathbb{M}_i}$.
Likewise, $\nodes{j}$ provides the set of nodes associated with stage $j$.
We refer to the node $m = 0$ at stage $j=0$ as the root node and the nodes at stage $j=J$ as leaf nodes.
The set of all non-root nodes is $\mathbb{L}_i = \mathbb{M}_i \setminus \{0\}$.
Each node $m \in \nodes{j}$ at stage $j\in\N_{[0, J-1]}$ is connected to a set of child nodes at stage $j+1$ which are accessible via $\child{m}$.
Vice versa, all nodes $m \in \nodes{j}$ at stage $j\in\N_{[1, J]}$ are reachable from one unique ancestor node at stage $j-1$, denoted $\anc{m}$.
The probability to visit node $m$ is $\pi^{(m)} \in (0, 1] \subset \mathbb{R}$.
By construction of the tree, we have that
\begin{subequations}
\begin{align}
	\smallSum\limits_{m \in \nodes{j}} \pi^{(m)} &= 1 \phantom{\pi^{(m)}} ~~ \forall j \in \mathbb{N}_{[0,J]} \quad \text{and} \\
	\smallSum\limits_{m_{+} \in \child{m}} \pi^{(m_+)} &= \pi^{(m)} \phantom{1} ~~ \forall m \in \mathbb{M}_i \setminus \nodes{J}.
\end{align}
\end{subequations}

\begin{figure}
	\vspace{0.1em}
	\includetikz{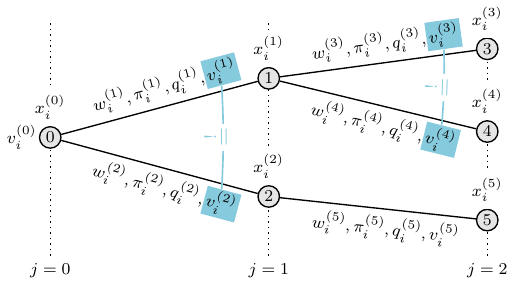}%
	\caption{Scenario tree based on \cite{HSB+2015} with nonanticipativity constraints.
	}
	\label{fig:scenarioTree}
\end{figure}

Each node $m \in \mathbb{M}_i$ is associated with a vector of states $x_i^{(m)}$ (see \cref{fig:scenarioTree}).
Moreover,
$v_i^{(m)}$ represents the control inputs,
$w_i^{(m)}$ the uncertain inputs and
$q_i^{(m)}$ the auxiliary variables
associated with node $m \in \mathbb{L}_i$.

\subsubsection{Nonanticipativity}

In the upcoming stochastic optimization problems, we aim to make control decisions $\ivar{v}{i}{(m)}$ at each node $m \in \mathbb{L}_i$.
In this context, it is important to move causally in time:
When making control decisions, we do not know exactly which uncertain input will occur.
We only know possible outcomes and associated probabilities.
Thus, we need make \emph{one} decision for \emph{all} uncertain inputs that originate from the same ancestor.
This must be accounted for, e.g., by the equality constraint
\begin{equation}\label{eq:nonanticipativity}
	\ivar{v}{i}{(m)} = \ivar{v}{i}{(n)} \quad \forall n \in \child{\anc{m}}.
\end{equation}
In \cref{fig:scenarioTree}, this so-called nonanticipativity constraint (see, e.g., \cite{DHB2018,SDR2014,RW1976}) is illustrated in light blue.
Here, at stage $j=1$, we can only make one control decision $\ivar{v}{i}{(1)} = \ivar{v}{i}{(2)}$ for the possible uncertain inputs $\ivar{w}{i}{(1)}$ and $\ivar{w}{i}{(2)}$.
At stage $j=2$, when finding control actions for nodes $3$ and $4$, we can only use the information that is available up to node $1$, i.e.,
we consider the case associated with uncertain input $\ivar{w}{i}{(1)}$.
Here, we need to make one decision $\ivar{v}{i}{(3)} = \ivar{v}{i}{(4)}$ for possible uncertain inputs $\ivar{w}{i}{(3)}$ and $\ivar{w}{i}{(4)}$.
Since node $5$ is a singleton, no restrictions of the form \eqref{eq:nonanticipativity} apply.

\begin{remark}
In problem formulations, nonanticipativity can be included in different ways.
One possibility is replacing $\ivar{v}{i}{(n)}$ by $\ivar{v}{i}{(m)}$ for all $n \in \child{\anc{m}}$ when formulating the problems.
Another way is to include \eqref{eq:nonanticipativity} as a regular constraint.
In this case, the pre-solve stage in off-the-shelf solvers, would often perform aforementioned replacement internally.
Either way, \eqref{eq:nonanticipativity} reduces the complexity of associated problems since less decision variables must be considered. \qed
\end{remark}

\subsection{Single \acl{mg} model}
\label{sec:mgModel}

In what follows, a model of a single \ac{ac} \ac{mg} is derived.
It is based on \cite{HSB+2015} and includes \acp{res}, storage and conventional units as well as a \ac{pcc}.

\subsubsection{\acsp{res}}
\begin{subequations}\label{eq:resBounds}
The forecasts of power and setpoint must lie within the bounds $p_{r,i}^{\min}\in \mathbb{R}_{\geq 0}^{R_i}$ and $p_{r,i}^{\max} \in \mathbb{R}_{> 0}^{R_i}$, i.e.,
\begin{align}
&p_{r,i}^{\min} \leq u_{r,i}^{(m)} \leq p_{r,i}^{\max}, \\
&p_{r,i}^{\min} \leq p_{r,i}^{(m)} \leq p_{r,i}^{\max},
\end{align}
for all $m \in \mathbb{L}_i$ \cite{HSB+2015}.
\end{subequations}
Additionally, the power of each \ac{res} can be limited by its setpoint.
The limit comes to bear if the weather-dependent available renewable power exceeds the setpoint.
If it lies below the setpoint, then the power follows the available infeed.
Using the element-wise $\min$ operator allows to describe this for all $m \in \mathbb{L}_i$ by \cite{HSB+2015}
\begin{equation}
p_{r,i}^{(m)}
= \min (
u_{r,i}^{(m)}
,
w_{r,i}^{(m)}
).
\label{eq:resBigM}
\end{equation}
Note that \eqref{eq:resBigM} can be easily transformed into a set of affine constraints using the auxiliary decision variable $\delta_{r,i}^{(m)}$ (see, e.g., \cite[Lemma 3.3.6]{Han2021}) which makes it appropriate for \acp{mip}.

\subsubsection{Storage units}
\label{sec:gridModel:storageUnits}

The forecasts of setpoints and power are limited by $p_{s,i}^{\min} \in \mathbb{R}_{<0}^{S_i}$ and $p_{s,i}^{\max} \in \mathbb{R}_{>0}^{S_i}$, i.e.,
\begin{subequations}
\label{eq:storageBounds}
\begin{align}
p_{s,i}^{\min} &\leq u_{s,i}^{(m)} \leq p_{s,i}^{\max}, \\
p_{s,i}^{\min} &\leq p_{s,i}^{(m)} \leq p_{s,i}^{\max},
\end{align}
\end{subequations}
for all $m \in \mathbb{L}_i$ \cite{HSB+2015}.
Moreover, the storage units exhibit dynamics which can be forecast for all $m \in \mathbb{L}_i$ via
\begin{equation}
	\ivar{x}{i}{(m)} = \ivar{x}{i}{(m_-)} - T_s p_{s,i}^{(m)},
	\label{eq:storageDynamics}
\end{equation}
with $m_- \in \anc{m}$ \cite{HSB+2015}.
Here, $x_i^{(0)} = x_{i}(k)$ is the measured state at the current discrete time instant $k\in\N_0$.
It is desired to keep $\ivar{x}{i}{(m)}$ above $\ivar{x}{i}{\min} \in \mathbb{R}_{\geq 0}^{S_i}$ and below $\ivar{x}{i}{\max} \in \mathbb{R}_{>0}^{S_i}$.
To ensure feasibility, $\ivar{\sigma}{i}{(m)} \in \mathbb{R}_{\geq 0}^{S_i}$ is used to form the soft constraints
\begin{equation}
x_{i}^{\min} - \ivar{\sigma}{i}{(m)}
	\leq x_{i}^{(m)} \leq
	x_i^{\max} + \ivar{\sigma}{i}{(m)}
\label{eq:energyBounds}
\end{equation}
for all $m \in \mathbb{L}_i$
which are completed by adding a penalty on nonzero values of $\ivar{\sigma}{i}{(m)}$ to the objective in \cref{sec:costFunction}.

\subsubsection{Conventional units}

\begin{subequations}
	\label{eq:conventionalBounds}
	The forecasts of power and setpoints must lie within
	$p_{t,i}^{\min} \in \mathbb{R}_{>0}^{T_i}$ and
	$p_{t,i}^{\max} \in \mathbb{R}_{>0}^{T_i}$
	for enabled units.
	With the Boolean input $\delta_{t,i}^{(m)}$ that indicates if units are enabled or disabled, this can be formulated as \cite{HSB+2015}
	\begin{align}
		\diag\left( p_{t,i}^{\min} \right) \delta_{t,i}^{(m)}
		\leq u_{t,i}^{(m)}
		\leq \diag\left( p_{t,i}^{\max} \right) \delta_{t,i}^{(m)}, \\
		\diag\left( p_{t,i}^{\min} \right) \delta_{t,i}^{(m)}
		\leq p_{t,i}^{(m)}
		\leq \diag\left( p_{t,i}^{\max} \right) \delta_{t,i}^{(m)}.
	\end{align}
\end{subequations}

\subsubsection{Power sharing between grid-forming units}
\label{sec:gridModel:powerSharing}

We assume that grid-forming storage and conventional units change their power in presence of fluctuations in a given proportional manner.
This is typically achieved via suitable low-level control schemes (see, e.g., \cite{SHK+2017}) and can be modeled by \cite{Han2021}
\begin{subequations}\label{eq:powerSharingInequalities}
	\begin{align}
		\label{eq:powerSharingConventional}
		K_{t,i} (\ivar{p}{t,i}{(m)} - \ivar{u}{t,i}{(m)} )
		&= \ivar{\rho}{i}{(m)} \ivar{\delta}{t,i}{(m)}, \\
		\label{eq:powerSharingStorage}
		K_{s,i} ( \ivar{p}{s,i}{(m)} - \ivar{u}{s,i}{(m)} )
		&= \ivar{\rho}{i}{(m)} \mathbf{1}_{S_i},
	\end{align}
\end{subequations}
for all $m \in \mathbb{L}_i$
with the auxiliary decision variable $\rho_i^{(m)}$ and the droop coefficient matrices
\begin{align*}
	K_{t,i} &= \diag( [ 1/\chi_{i,1} \quad \cdots \quad 1/\chi_{i,T_i} ]\tran ),\\
	K_{s,i} &= \diag( [ 1/\chi_{i,(T_i+1)} \quad \cdots \quad 1/\chi_{i,(T_i+S_i)} ]\tran ).
\end{align*}
Here, the design parameters $\chi_{i,l} \in \mathbb{R}_{>0}$, $l\in\N_{[1,T_i+S_i]}$ can be chosen for example according to the units' nominal power.

The term $\ivar{\rho}{i}{(m)} \ivar{\delta}{t,i}{(m)}$ in \eqref{eq:powerSharingConventional} models that only enabled units can participate in power sharing.
It can be easily transformed into a set of affine constraints (see, e.g., \cite[Lemma 3.3.5]{Han2021}) which renders it appropriate for mixed-integer formulations.

\subsubsection{\acs{pcc} power}
\label{sec:mg:pcc}

We assume that the low-level control schemes keep the power at the \ac{pcc} at a desired value $p_{g,i}(j)$ in presence of uncertain demand and renewable infeed.
Thus, all fluctuations are covered inside each \ac{mg} $i$ and power exchange over the utility grid is \emph{not} affected by uncertainties, i.e., it is certain.
In the model, this is captured by only considering \emph{one} predicted \ac{pcc} power value $p_{g,i}(j)$ which is identical for all nodes at stage $j \in \mathbb{N}_{[1,J]}$.
This power is limited by $p_{g,i}^{\min} \in \mathbb{R}_{\leq 0}$ and $p_{g,i}^{\max} \in \mathbb{R}_{\geq 0}$, i.e.,
\begin{equation}
	\label{eq:powerExchangeBounds}
	p_{g,i}^{\min} \leq p_{g,i}(j) \leq p_g^{\max}.
\end{equation}

\subsubsection{Local power balance}
Within each \ac{mg} $i \in \mathbb{I}$, a local power balance holds.
For all $m \in \mathbb{L}_i$, it is modeled by
\begin{multline}
0 = \mathbf{1}_{R_i}\tran p_{r,i}^{(m)}
+ \mathbf{1}_{T_i}\tran p_{t,i}^{(m)}
+ \mathbf{1}_{S_i}\tran p_{s,i}^{(m)}
+ \mathbf{1}_{D_i}\tran w_{d,i}^{(m)} \\
+ p_{g,i} (\stage{m}).
\label{eq:localPowerBalance}
\end{multline}

\subsection{Grid model}
\label{sec:utilityModel}

\begin{subequations}\label{eq:transmissionConstraints}
Let us collect the \ac{pcc} power forecast $p_{g,i}(j)$ of all \acp{mg} at stage $j \in \mathbb{N}_{[1,J]}$ in $p_{g}(j) = [ p_{g,1}(j) ~ \cdots ~ p_{g,I}(j) ]\tran$.
Negative values of $p_{g,i}(j)$ indicate that \ac{mg} $i$ injects power into the grid, whereas a positive values indicates a consumption of power.
Let us further collect the power flow over ${E\in\N}$ transmission lines in $p_{e}(k) = [ p_{e,1}(k) ~ \cdots ~ p_{e,E}(k) ]\tran$.
Assuming inductive short to medium length power lines allows us to employ the \ac{dc} power flow approximations for \ac{ac} grids \cite{PMVDB2005}.
Motivated by \cite{HBR+2018}, we describe the power flow over the lines by
\begin{equation}
	p_{e}(j) = \tilde{F} \cdot p_g(j).
	\label{eq:powerFlow}
\end{equation}
Moreover, a global power equilibrium of the form
\begin{equation}
	0 = \mathbf{1}_I\tran p_g (j)
\label{eq:globalPowerBalance}
\end{equation}
must hold.
The power that can be transmitted via the grid is limited by $p_{e}^{\min} \in \mathbb{R}_{<0}^{E}$ and $p_{e}^{\max} \in \mathbb{R}_{>0}^{E}$, i.e.,
\begin{equation}
	p_{e}^{\min} \leq p_{e}(j) \leq p_{e}^{\max}.
	\label{eq:powerFlowBounds}
\end{equation}
\end{subequations}

%% file: problemFormulation.tex

\section{Problem formulation}
\label{sec:problemFormulation}

Based on the model from the previous chapter, we can now formulate \ac{mpc} problems for the operation of interconnected \acp{mg}.
We start by discussing the operation costs.

\subsection{Costs of individual \acsp{mg}}
\label{sec:costFunction}

In what follows, we define the units' costs for all nodes $m\in\mathbb{L}_i$ in the scenario trees of each \ac{mg} $i\in\mathbb{I}$.

\subsubsection{\acsp{res}}

The desire for high renewable infeed is considered by penalizing deviations from the rated power $p_{r,i}^{\max}$.
With weight $c_{r,i} \in \mathbb{R}_{>0}^{R_i}$ the cost is formulated as \cite{HBR+2018}
\begin{equation}
\ell_{r,i}^{(m)} = \norm{\diag(c_{r,i})(\ivar{p}{r,i}{\max} - \ivar{p}{r,i}{(m)})}_2^2.
\end{equation}

\subsubsection{Storage units}

The costs of the storage units are
\begin{equation}
\ivar{\ell}{s,i}{(m)} =
	 \norm{ \diag(c_{s,i})  \ivar{p}{s,i}{(m)} }_2^{2}
	 + \norm{ \diag(c_{\sigma,i}) \ivar{\sigma}{i}{(m)} }_2^2
\end{equation}
with $c_{s,i} \in \mathbb{R}_{>0}^{S_i}$ and $c_{\sigma,i} \in \mathbb{R}_{>0}^{S_i}$.
The first term reflects conversion losses,
the second term is a penalty for nonzero values of $\ivar{\sigma}{i}{(m)}$, i.e., energy values above ${x}_i^{\max}$ or below ${x}_{i}^{\min}$ (see \cref{sec:gridModel:storageUnits}).

\subsubsection{Conventional generators}

Following \cite{JMK2012}, the operating costs of the conventional units are modeled using a quadratic function.
With $c_{t,i},\,c_{t,i}^{\prime},\,c_{t,i}^{\prime\prime} \in \mathbb{R}_{>0}^{T_i}$, it reads
\begin{equation}
\ell_{t,i}^{(m)}
	= c_{t,i}\tran\delta_{t,i}^{(m_-)}
	+ c_{t,i}^{\prime \mathstrut\scriptscriptstyle\top}p_{t,i}^{(m)}
	+ \norm{\diag(c_{t,i}^{\prime\prime})p_{t,i}^{(m)}}_2^2.
\end{equation}

Switching actions incur maintenance costs.
This is accounted for by a term that is nonzero if units are enabled or disabled from node ${m_{-} = \anc{m}}$ to node $m$, i.e.,
\begin{equation}
\ell_{sw,i}^{(m)} = \norm{\diag(c_{sw,i})(\delta_{t,i}^{(m_-)} - \delta_{t,i}^{(m)})}_{2}^{2},
\end{equation}
with $c_{sw,i} \in \mathbb{R}_{>0}^{T_i}$ \cite{HSB+2015}.
Here, $\delta_{t,i}^{(0)}$ equals the measured $\delta_{t,i}(k)$ at the current time instant $k$.

\subsubsection{Control effort of lower layers}
Large differences between power and setpoint can lead to additional control effort at the lower layers.
Therefore, it is desirable to keep the power close to the setpoints.
Since $\ivar{\rho}{i}{(m)}$ is correlated with deviations of the power from the setpoints, this desire can be reflected by a penalty, with weight $c_{\rho,i} \in \mathbb{R}_{>0}$, of the form
\begin{equation}
\ivar{\ell}{\rho,i}{(m)} = c_{\rho,i} \, ( \ivar{\rho}{i}{(m)} )^2.
\end{equation}

\subsubsection{Trading cost}

All fluctuations are assumed to be covered locally inside each \ac{mg}.
As a result, $p_g$ and consequently the cost for power exchange, i.e.,
\begin{equation}\label{eq:tradingCost}
	\ivar{\ell}{g,i}{(m)} =  c_{g,i} p_{g,i}(\stage{m})
		+ \ivar{c}{g,i}{\prime} |p_{g,i}(\stage{m})|
\end{equation}
with $\ivar{c}{g,i}{\prime}$, $\ivar{c}{g,i}{\prime\prime} \in \mathbb{R}_{>0}$ \cite{HBR+2018}, are certain.
The first part of \eqref{eq:tradingCost} represents a given price for energy and the second part a fixed cost per absolute value of traded energy.

\subsubsection{Single \acs{mg} cost}

By construction of the scenario tree (see \cref{sec:scenario-structure}), the expected cost is given by the sum over all nodes $m\in\mathbb{L}_i$, weighted with probability $\ivar{\pi}{i}{(m)}$, i.e.,
\begin{multline}
\ell_{i} = \smallSum\limits_{m \in \mathbb{L}_i} \ivar{\pi}{i}{(m)} \big(
	\ivar{\ell}{r,i}{(m)} + \ivar{\ell}{s,i}{(m)} + \ivar{\ell}{t,i}{(m)} + \ivar{\ell}{sw,i}{(m)}
	\big. \\
	\big. + \ivar{\ell}{\rho,i}{(m)} + \ivar{\ell}{g,i}{(m)} \big) \gamma^{\stage{m}}. \label{eq:expectedCost}
\end{multline}
Here, $\gamma \in (0,1)$ is used to emphasize near future decisions.
Note that \cref{eq:expectedCost} equals the sum of expected stage costs over all $j\in\N_{[1, J]}$ in the scenario tree (see, e.g., \cite[Ch.~10]{Han2021}).

\subsection{Cost for power transmission}

In addition to \eqref{eq:expectedCost}, costs for transmitting power are considered.
Assuming that they increase with line losses, i.e., quadratically with the transmitted power $p_{e}$,
allows us to deduce the transmission costs over prediction horizon $J$ with diagonal matrix $C_{e}\in\R^{E \times E}_{\geq 0}$ as \cite{HBR+2018}
\begin{equation}
	\ell_{e} = \smallSum\limits_{j=1}^{J} p_{e}\tran(j) C_{e} p_{e}(j) \cdot \gamma^{j}.
\end{equation}

\subsection{\acs{mpc} problem formulation}
\label{sec:centralProblem}

For each \ac{mg} $i\in\mathbb{I}$, we collect
all control inputs in $\mathbf{V}_i = [\ivar{v}{i}{(1)} ~ \cdots ~ \ivar{v}{i}{(M_i-1)} ]$,
states in $\mathbf{X}_i = [\ivar{x}{i}{(0)} ~ \cdots ~ \ivar{x}{i}{(M_i-1)}]$
and
auxiliary variables in $\mathbf{Q}_i = [\ivar{q}{i}{(1)} ~ \cdots ~ \ivar{q}{i}{(M_i-1)}]$.
Similarly, we form $\mathbf{P}_{e} = [ p_{e}(1) ~ \cdots ~ p_{e}(J) ]$ and $\mathbf{P}_g = [ p_g(1) ~ \cdots ~ p_g(J) ]$.
Using these decision variables, we can formulate the following stochastic optimization problem.

\begin{problem}[Central mixed-integer problem]
\label{prob:mip}
\begin{equation*}
	\minimize_{\substack{
		\mathbf{V}_1,\ldots,\mathbf{V}_I \\
		\mathbf{X}_1,\ldots,\mathbf{X}_I \\
		\mathbf{Q}_1,\ldots,\mathbf{Q}_I \\
		\mathbf{P}_{e}, \mathbf{P}_g}}
		\ell_{e} + \smallSum\limits_{i\in\mathbb{I}} \ell_{i}
\end{equation*}
subject to

	\cref{eq:nonanticipativity,eq:storageDynamics,eq:resBounds,eq:resBigM,eq:energyBounds,eq:storageBounds,eq:conventionalBounds,eq:powerSharingInequalities,eq:powerExchangeBounds,eq:localPowerBalance,eq:transmissionConstraints} for all $m\in\mathbb{L}_i$ \\ \indent
	as well as initial conditions $\ivar{x}{i}{(0)} = x_{i}(k)$, $\ivar{\delta}{t,i}{(0)} =  \delta_{t,i}(k)$ \\ \indent
	for all $i \in \mathbb{I}$. \qed
\end{problem}

\begin{remark}[Robustness]
In \cref{prob:mip}, power and energy, i.e., $\mathbf{X}_i$ and $\mathbf{Q}_i$, are monotone in the uncertain input \cite[Remark 9.2.4]{Han2021}.
Therefore, the formulation is stage-wise robust to uncertain inputs which are in the convex hull formed by $w_i^{(m)}$, $m\in\nodes{j}$.
In \cref{sec:caseStudy}, this property is exploited by adding extreme scenarios with low probabilities to the tree, which extends the convex hull and allow to increase robustness while keeping good performance when optimizing the expected cost. \qed
\end{remark}

The decision variables $\delta_{t,i} \in \{0,1\}^{T_i}$ and $\delta_{r,i}\in \{0,1\}^{R_i}$ are Boolean, which renders \cref{prob:mip} a \acl{miqp}.
In order to find a solution along the lines of \cite{HBR+2018}, we relax $\delta_{t,i}^{(m)} \in [0,1]^{T_i}$ and $\delta_{r,i}^{(m)} \in [0,1]^{R_i}$, which allows to deduce the following quadratic formulation.

\begin{problem}[Central relaxed problem]
\label{prob:relaxed}
\begin{equation*}
	\minimize_{\substack{
		\mathbf{V}_1,\ldots,\mathbf{V}_I \\
		\mathbf{X}_1,\ldots,\mathbf{X}_I \\
		\mathbf{Q}_1,\ldots,\mathbf{Q}_I \\
		\mathbf{P}_{e}, \mathbf{P}_g}}
	\ell_{e} + \smallSum\limits_{i\in\mathbb{I}} \ell_{i}
\end{equation*}
subject to

	\cref{eq:nonanticipativity,eq:storageDynamics,eq:resBounds,eq:resBigM,eq:energyBounds,eq:storageBounds,eq:conventionalBounds,eq:powerSharingInequalities,eq:powerExchangeBounds,eq:localPowerBalance,eq:transmissionConstraints} for all $m\in\mathbb{L}_i$ \\ \indent
	as well as initial conditions $\ivar{x}{i}{(0)} = x_{i}(k)$, $\ivar{\delta}{t,i}{(0)} =  \delta_{t,i}(k)$ \\ \indent
	with $\delta_{t,i}^{(m)} \in [0,1]^{T_i}$ and $\delta_{r,i}^{(m)} \in [0,1]^{R_i}$ for all $m\in\mathbb{L}_i$ \\ \indent
	and all $i \in \mathbb{I}$. \qed
\end{problem}

%% file: distributedSolution.tex

\section{Distributed solution}
\label{sec:distributedSolution}

\begin{figure}
	\vspace{0.5em}
  \centering
  \input{figures/communicationStructure}%
  \caption{Underlying communication structure in \cref{algo:distributed}.}
  \label{fig:communicationStructure}
\end{figure}
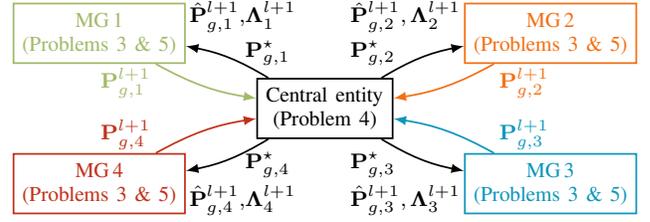

In this section, we describe how to find a (not necessarily optimal) distributed solution to \cref{prob:mip} using \cref{algo:distributed}.
In detail, we first find a solution to \cref{prob:relaxed} by alternately solving a subproblem at a central coordinator and local subproblems at the \acp{mg} (see \cref{fig:communicationStructure}).
In a second step, local \acp{mip} are solved at each \acp{mg} for the given \ac{pcc} power from the previous step.

Consider a copy of $\mathbf{P}_g$ denoted $\mathbf{\hat{P}}_g = [ \hat{p}_g(1) ~ \cdots ~ \hat{p}_g(J) ]$.
The decision variables in $\mathbf{P}_g$ are used at the local \acp{mg} and the ones in $\mathbf{\hat{P}}_g$ at the central entity.
Both have to be found such that
\begin{equation}
\label{eq:pg-equality}
\mathbf{P}_g - \mathbf{\hat{P}}_g = \mathbf{0}_{I\times J} \; .
\end{equation}

Let us form the vector of Lagrange multipliers of all \acp{mg} at prediction step $j$ as $\lambda(j) = \left[ \lambda_1(j) ~ \cdots ~ \lambda_I(j) \right]\tran \in \R^I$
and collect them in $\boldsymbol{\Lambda} = \left[ \lambda(1) ~ \cdots ~ \lambda(J) \right]$.
We refer
to row $i$ of $\boldsymbol{\Lambda}$ as $\boldsymbol{\Lambda}_{i}$,
to row $i$ of $\mathbf{P}_{g}$ as $\mathbf{P}_{g,i}$ and
to row $i$ of $\mathbf{\hat{P}}_{g}$ as $\mathbf{\hat{P}}_{g,i}$.
This allows us to formulate the augmented Lagrangian of \cref{prob:relaxed} with fixed parameter $\kappa \in \mathbb{R}_{>0}$,
\begin{multline}
	\mathcal{L}_{\kappa}(\mathbf{P}_{g},\mathbf{\hat{P}}_{g},\boldsymbol{\Lambda})
	=
	\ell_{e} + \smallSum\limits_{i\in\mathbb{I}}
	\big(
	\ell_{i} +
	\boldsymbol{\Lambda}_{i} (\mathbf{P}_{g,i} - \mathbf{\hat{P}}_{g,i})\tran \\ +
	\nicefrac{\kappa}{2} \| \mathbf{P}_{g,i} - \mathbf{\hat{P}}_{g,i} \|_{2}^2
	\big).
\end{multline}

Similar to \cite{HBR+2018}, we intend to find a distributed solution to \cref{prob:relaxed} based on the augmented Lagrangian.
Therefore, we pose the following optimization problems.

\begin{problem}[Local \acs{admm} problem at \acs{mg} $i\in\mathbb{L}$]
\label{prob:distributedLocal}
\[
	\ivar{\mathbf{P}}{g,i}{l+1} \in \argmin_{\mathbf{V}_i, \mathbf{X}_i, \mathbf{Q}_i, \mathbf{P}_{g,i}}
		\ell_i +
		{\boldsymbol{\Lambda}}_{i}^{l} \mathbf{P}_{g,i}\tran
		+ \nicefrac{\kappa}{2} \| \mathbf{P}_{g,i} - \ivar{\mathbf{\hat{P}}}{g,i}{l} \|_2^2
\]
subject to

	\cref{eq:nonanticipativity,eq:storageDynamics,eq:resBounds,eq:resBigM,eq:energyBounds,eq:storageBounds,eq:conventionalBounds,eq:powerSharingInequalities,eq:powerExchangeBounds,eq:localPowerBalance} for all $m\in\mathbb{L}_i$ \\ \indent
	as well as initial conditions $\ivar{x}{i}{(0)} = x_{i}(k)$, $\ivar{\delta}{t,i}{(0)} =  \delta_{t,i}(k)$ \\ \indent
	with $\delta_{t,i}^{(m)} \in [0,1]^{T_i}$, $\delta_{r,i}^{(m)} \in [0,1]^{R_i}$ for all $m\in\mathbb{L}_i$. \qed
\end{problem}

\vspace{0.5\baselineskip}

\begin{problem}[Central \acs{admm} problem at coordinator]
\begin{samepage}
\label{prob:distributedCentral}
\[
	\ivar{\hat{\mathbf{P}}}{g}{l+1} \in \argmin_{\hat{\mathbf{P}}_{g},\mathbf{P}_{e}} ~
	\ell_{e} - \smallSum\limits_{i\in\mathbb{I}}
	\big(
	\boldsymbol{\Lambda}_{i}^{l} \mathbf{\hat{P}}_{g,i}\tran
	+ \nicefrac{\kappa}{2} \| \ivar{\mathbf{P}}{g,i}{l+1} - \mathbf{\hat{P}}_{g,i} \|_2^2
	\big)
\]
subject to  \\ \indent
	\cref{eq:powerExchangeBounds,eq:transmissionConstraints} using $\hat{p}_g(j)$ instead of ${p}_g(j)$ \\ \indent
	for all $j \in \mathbb{N}_{[1,J]}$. \qed
\end{samepage}
\end{problem}

By alternately solving \cref{prob:distributedCentral,prob:distributedLocal} and updating Lagrange multipliers (see \eqref{eq:lagrangeMultipliers} in \cref{algo:distributed}) we can find an optimal solution to \cref{prob:relaxed} using the \ac{admm}.
The resulting optimal values $\ivar{\mathbf{P}}{g,i}{\star}$ are then used to find feasible (but not necessarily optimal) solutions to \cref{prob:mip} via local mixed-integer updates at all \acp{mg}.
In these updates,
$\mathbf{P}_g$ is fixed to $\ivar{\mathbf{P}}{g,i}{\star}$ and
Boolean values for $\delta_{t,i}$ and $\delta_{r,i}$
are again considered.
The associated problem reads as follows.

\begin{problem}[Mixed-integer update at \acs{mg} $i\in\mathbb{L}$]
\label{prob:mipUpdate}
\begin{equation*}
	\minimize_{
		\mathbf{V}_i,
		\mathbf{X}_i,
		\mathbf{Q}_i
	} ~ \ell_i
\end{equation*}
subject to

	\cref{eq:nonanticipativity,eq:storageDynamics,eq:resBounds,eq:resBigM,eq:energyBounds,eq:storageBounds,eq:conventionalBounds,eq:powerSharingInequalities,eq:localPowerBalance} for all $m\in\mathbb{L}_i$, \\ \indent
	as well as initial conditions $\ivar{x}{i}{(0)} = x_{i}(k)$, $\ivar{\delta}{t,i}{(0)} =  \delta_{t,i}(k)$ \\ \indent
	with fixed $\mathbf{P}_{g,i} = \ivar{\mathbf{P}}{g,i}{\star}$. \qed
\end{problem}

\Cref{algo:distributed} allows us to find a hierarchical distributed solution using \cref{prob:distributedLocal,prob:distributedCentral,prob:mipUpdate} (see also \cref{fig:communicationStructure}).
Similar to \cite{HBR+2018}, a termination criterion that checks if
the change of Lagrange multipliers,
\ac{pcc} power,
and residuals
are all below a small $\epsilon\in\R_{>0}$.
If this termination criterion is not met at $l = l_{\max}$, then $\ivar{\hat{\mathbf{P}}}{g,i}{l_{\max}+1}$, which represents a feasible solution to \cref{prob:relaxed}, is used subsequently.
Finally, in step 3, local \acp{mip} are employed to find a feasible solution to \cref{prob:mip}.

\begin{remark}[Privacy]
	By design of \cref{algo:distributed}, the local \ac{mg} controllers only share the \ac{pcc} power forecast $\ivar{\mathbf{P}}{g,i}{l+1}$ with the central entity.
	Thus, no explicit information about the structure of the \ac{mg} in the form of constraints or cost function is shared with others which helps to preserve the privacy of the local \ac{mg} controllers.
	One can, however, think of ways to reconstruct parts of \cref{prob:distributedLocal} using malicious vectors $\ivar{\hat{\mathbf{P}}}{g,i}{l+1}$, $\boldsymbol{\Lambda}_{i}^{l+1}$.
	Finding formulations which are secure against such attacks is subject to future work. \qed 
\end{remark}

\begin{remark}[Suboptimality]
	The \ac{admm} part of \cref{algo:distributed} provides an optimal solution to \cref{prob:relaxed}.
	However, the overall scheme, composed of \ac{admm} and \cref{prob:mipUpdate} is not guaranteed to find an optimal solution to \cref{prob:mip}.
	That being said, in the simulations performed in \cref{sec:caseStudy}, the algorithm was found to provide results with a small suboptimality gap to the original \cref{prob:mip}. \qed
\end{remark}

\begin{remark}[Feasibility]
	Let us assume that each \ac{mg} is configured such that storage and conventional units together can always serve all possible load values present in the forecast scenarios.
	Thus, an islanded operation where all \acp{res} are set to provide zero power always represents a feasible solution to \cref{prob:mipUpdate}.
	This allows us to deal with \ac{admm} solutions $\ivar{\mathbf{P}}{g,i}{\star}$ which lead to infeasible formulations of \cref{prob:mipUpdate} at \ac{mg} $i$ by re-executing \cref{algo:distributed} with zero \ac{pcc} power using $p_{g,i}^{\min} = p_{g,i}^{\max} = 0$.
	This results in $\ivar{\mathbf{P}}{g,i}{\star} = \mathbf{0}_{1\times J}$ which in turn lead to feasible formulations of \cref{prob:mipUpdate}.
	Another alternative is to employ schemes that ensure feasibility, such as the one presented in \cite{SHHR2019}.
	For a thorough discussion of feasibility in a related problem, the reader is kindly referred to \cite[Sec.~V.C]{HBR+2018}.
	 \qed
\end{remark}

\begin{algorithm}[h]
\caption{Hierarchical distributed algorithm}
\label{algo:distributed}
	\begin{enumerate}[label=\textbf{\arabic*.}, leftmargin=1.3em]
	\item \textbf{Initialize:} At time $k$, $\forall i \in \mathbb{I}$, measure $x_{i}(k)$, $\delta_{t,i}(k)$ and obtain scenario tree.
	\item \textbf{\acs{admm} loop:} \\
		\textbf{for} $l = 0, \ldots, l_{\max} \in \mathbb{N}$:
		\begin{enumerate}[label=(\roman*), leftmargin=1.75em]
		\item For all \ac{mg} $i \in \mathbb{I}$ (in parallel):
		\begin{itemize}[leftmargin=1.0em]
			\item Solve \cref{prob:distributedLocal} in parallel to obtain $\ivar{\mathbf{P}}{g,i}{l+1}$.
			\item Send $\ivar{\mathbf{P}}{g,i}{l+1}$ to central entity.
		\end{itemize}
		\item Central entity:
		\begin{itemize}[leftmargin=1.0em]
			\item Solve \cref{prob:distributedCentral} to obtain $\ivar{\hat{\mathbf{P}}}{g}{l+1}$.
			\item Update Lagrange multipliers:
			\begin{equation}
				\boldsymbol{\Lambda}^{l+1} = \boldsymbol{\Lambda}^{l} + \kappa \big(
					\ivar{\mathbf{P}}{g}{l+1} - \ivar{\hat{\mathbf{P}}}{g}{l+1}
					\big).
				\label{eq:lagrangeMultipliers}
			\end{equation}
			\item Communicate $\ivar{\hat{\mathbf{P}}}{g,i}{l+1}$ and $\boldsymbol{\Lambda}_{i}^{l+1}$ to all \ac{mg} $i \in \mathbb{I}$.
			\item Check termination criterion: \newline
				\textbf{if}
				$\big(| \boldsymbol{\Lambda}^{l} - \boldsymbol{\Lambda}^{l+1} | < \epsilon\big.$
				{and}
				$| \ivar{\mathbf{P}}{g,i}{l} - \ivar{\mathbf{P}}{g,i}{l+1}| < \epsilon$
				{and} \newline \phantom{\textbf{if(}}
				$\big. |\ivar{\mathbf{P}}{g,i}{l+1} - \ivar{\mathbf{\hat{P}}}{g,i}{l+1}| < \epsilon\big)$ {or}  $l = l_{\max}$, \newline
				\textbf{then} set $\ivar{\mathbf{P}}{g,i}{\star} = \ivar{\hat{\mathbf{P}}}{g,i}{l+1}$ and go to 3.
		\end{itemize}
		\end{enumerate}
	\item \textbf{Mixed-integer update:} For all \ac{mg} $i \in \mathbb{I}$ (in parallel):
	\begin{itemize}[leftmargin=1.0em]
		\item Solve \cref{prob:mipUpdate}.
	\end{itemize}
	\end{enumerate}
\end{algorithm}

\Cref{algo:distributed} represents the core of an \ac{mpc} scheme:
At each time instant $k$, new scenario trees and new measurements $x_{i}(k)$, $\delta_{t,i}(k)$ are obtained at each \ac{mg} $i\in\mathbb{I}$.
Then, \Cref{algo:distributed} is used to find control actions $v_i^{(1)}$ which are applied to the each \ac{mg} $i\in\mathbb{I}$.
At time $k+1$, new measurements and new scenario trees are obtained and the scheme is repeated in a receding horizon manner.

%% file: figures/communicationStructure.tex

\tikzset{external/export next=false}

\begin{tikzpicture}[font=\footnotesize,
line width = 0.65pt,
>=latex,
]

\node[black, draw, align=center] (centralEntity) at (0, 0)
	{Central entity \\ (\cref{prob:distributedCentral})};

\node[vgGreen, draw, align=center] (mg1) at (-3, 1)
	{\acs{mg}\,1 \\ (Problems \ref{prob:distributedLocal} \& \ref{prob:mipUpdate})};

\node[vgOrange, draw, align=center] (mg2) at (3, 1)
	{\acs{mg}\,2 \\ (Problems \ref{prob:distributedLocal} \& \ref{prob:mipUpdate})};

\node[vgLightBlue, draw, align=center] (mg3) at (3, -1)
	{\acs{mg}\,3 \\ (Problems \ref{prob:distributedLocal} \& \ref{prob:mipUpdate})};

\node[vgRed, draw, align=center] (mg4) at (-3, -1)
	{\acs{mg}\,4 \\ (Problems \ref{prob:distributedLocal} \& \ref{prob:mipUpdate})};

\node[black, anchor=west, inner sep=1pt] at (mg1.0) {$\begin{aligned} 
		\ivar{\hat{\mathbf{P}}}{g,1}{l+1},
		& \boldsymbol{\Lambda}_{1}^{l+1} \\ 
		& \ivar{\mathbf{P}}{g,1}{\star}
	\end{aligned}$};
\node[black, anchor=east, inner sep=1pt] at (mg2.180) {$\begin{aligned} 
		& \ivar{\hat{\mathbf{P}}}{g,2}{l+1},
		\boldsymbol{\Lambda}_{2}^{l+1} \\
		& \ivar{\mathbf{P}}{g,2}{\star} 
	\end{aligned}$};
\node[black, anchor=east, inner sep=1pt] at (mg3.180) {$\begin{aligned} 
		& \ivar{\mathbf{P}}{g,3}{\star} \\ 
		& \ivar{\hat{\mathbf{P}}}{g,3}{l+1},
		\boldsymbol{\Lambda}_{3}^{l+1}
	\end{aligned}$};
\node[black, anchor=west, inner sep=1pt] at (mg4.0) {$\begin{aligned} 
		& \ivar{\mathbf{P}}{g,4}{\star} \\ 
		\ivar{\hat{\mathbf{P}}}{g,4}{l+1},
		& \boldsymbol{\Lambda}_{4}^{l+1}
	\end{aligned}$};

\node[vgGreen, anchor=north, inner sep=1pt, xshift=10] at (mg1.270) {$\ivar{\mathbf{P}}{g,1}{l+1}$};
\node[vgOrange, anchor=north, inner sep=1pt, xshift=-10] at (mg2.270) {$\ivar{\mathbf{P}}{g,2}{l+1}$};
\node[vgLightBlue, anchor=south, inner sep=1pt, xshift=-10] at (mg3.90) {$\ivar{\mathbf{P}}{g,3}{l+1}$};
\node[vgRed, anchor=south, inner sep=1pt, xshift=10] at (mg4.90) {$\ivar{\mathbf{P}}{g,4}{l+1}$};

\draw[->, bend right=10] (centralEntity) to (mg1);
\draw[->, bend left=10] (centralEntity) to (mg2);
\draw[->, bend right=10] (centralEntity) to	(mg3);
\draw[->, bend left=10] (centralEntity) to (mg4);

\draw[vgGreen, ->, bend right=10] (mg1) to (centralEntity);
\draw[vgOrange, ->, bend left=10] (mg2) to (centralEntity);
\draw[vgLightBlue, ->, bend right=10] (mg3) to (centralEntity);
\draw[vgRed, ->, bend left=10] (mg4) to (centralEntity);

\end{tikzpicture}

%% file: caseStudy.tex

\section{Case study}
\label{sec:caseStudy}

The following study compares the closed-loop results obtained with \cref{algo:distributed} with those of the certainty equivalence approach from \cite{HBR+2018}.
We start with the simulation setup.

\subsection{Simulation setup}
\label{sec:simulationSetup}

\begin{table}
\vspace{0.8em}
\centering
\caption{Simulation parameters ($i = 1,\ldots,4$)}
\label{tab:parameters}
\begin{tabular}{c l c l}
\toprule
Parameter & Value & Parameter & Value\\
\cmidrule(lr){1-2} \cmidrule(lr){3-4}
$[ \ivar{{x}}{i}{\min}, ~\ivar{{x}}{i}{\max}] $ & $[0.2,\, 6]\,\unit{puh}$ &
	$[ \ivar{p}{\mathrm{e},i}{\min}, ~ \ivar{p}{\mathrm{e},i}{\max} ] $ & $[-1,\, 1]\,\unit{pu} $ \\[4pt]

$[ \ivar{p}{t,i}{\min}, ~ \ivar{p}{t,i}{\max} ] $ & $[0.4,\, 1]\,\unit{pu}$ &
	$[ \ivar{p}{g,i}{\min},~ \ivar{p}{g,i}{\max} ] $ & $[-1,\, 1]\,\unit{pu}$ \\[4pt]

$[ \ivar{p}{s,i}{\min}, ~ \ivar{p}{s,i}{\max}] $ & $[-1,\, 1]\,\unit{pu}$ &
	$K_{s,i} = K_{t,i}$ & $1$\\[4pt]

$[ \ivar{p}{r,i}{\min}, ~ \ivar{p}{r,i}{\max}]$ & $[0,\, 2]\,\unit{pu}$ &
$[{\delta}_{1}{(0)}, \ldots, {\delta}_{4}{(0)}]$ &
	$[0,\, 0,\, 0,\, 0]$ \\
\bottomrule
\end{tabular}
\end{table}

In what follows, the network in \cref{fig:interconnectedMgScheme} with the simulation parameters in \cref{tab:parameters} is considered.
As initial state, $[{x}_{1} (0), \ldots, {x}_{1}{(0)}] = [1, 3.4, 2.9, 5.6]\,\unit{puh}$ and as line admittances $y_{e, i} = \unit[20]{pu}$, $i\in\N_{[1,4|}$ were used.
Additionally, the weights in \cref{tab:weights} were employed in the cost functions.

The simulation was executed for 336 steps with a sampling interval of $\unit[30]{min.}$, resulting in a total duration of 7 days.
In the \ac{mpc} schemes, a forecast prediction horizon of $12$ sampling intervals was considered.
The four scenario trees, i.e., one for each \ac{mg} to formulate \cref{prob:distributedLocal,prob:mipUpdate}, were generated from $500$ independent forecast scenarios at each execution of the controller.
These scenarios were deduced from Monte-Carlo simulations assuming wind speed, irradiance and demand forecasts with normally distributed residuals. 
A branching factor of $6$ at stage $1$, $2$ at stage $2$ and $1$ for all subsequent stages was considered.
Each tree was generated using forward selection \cite{HR2003}.
To each tree the largest and smallest forecast (at the first prediction step) from the collections of independent scenarios were added to increase robustness (see \cite[Section 12.1.2]{Han2021}).

\begin{table}
\centering
\caption{Weights in cost functions of all \acsp{mg} $(i = 1,\ldots,4)$.}
\label{tab:weights}
\begin{tabular}{c l c l}
	\toprule
	Weight & Value & Weight & Value\\
	\cmidrule(lr){1-2} \cmidrule(lr){3-4}
	$c_{r,i}$ & 1 $\nicefrac{1}{\textrm{pu}}$ & $c_{sw,i}$ & 0.1\\
	$c_{s,i}$ & 0.2236 $\nicefrac{1}{\textrm{pu}}$ & $c_{g,i}$ & 0.5 $\nicefrac{1}{\textrm{pu}}$\\
	$c_{t,i}$ & 0.1178 & $c_{g,i}^{\prime}$ & 0.1 $\nicefrac{1}{\textrm{pu}}$\\
	$c_{t,i}^{\prime}$ & 0.751 $\nicefrac{1}{\textrm{pu}}$ & $C_{\mathrm{e}}$ & $0.1 \cdot\operatorname{diag}([1,2,3,6]^\top) \nicefrac{1}{\textrm{pu}^2}$\\
	$c_{t,i}^{\prime\prime}$ & 0.0693 $\nicefrac{1}{\textrm{pu}}$ & $c_{\sigma,i}$ & 1000\\
	$c_{\rho,i}$ & 0.05 & $\gamma$ & 0.95\\
	\bottomrule
\end{tabular}
\end{table}

\begin{figure*}[t]
  \centering
	\vspace{0.5em}
  \includetikz{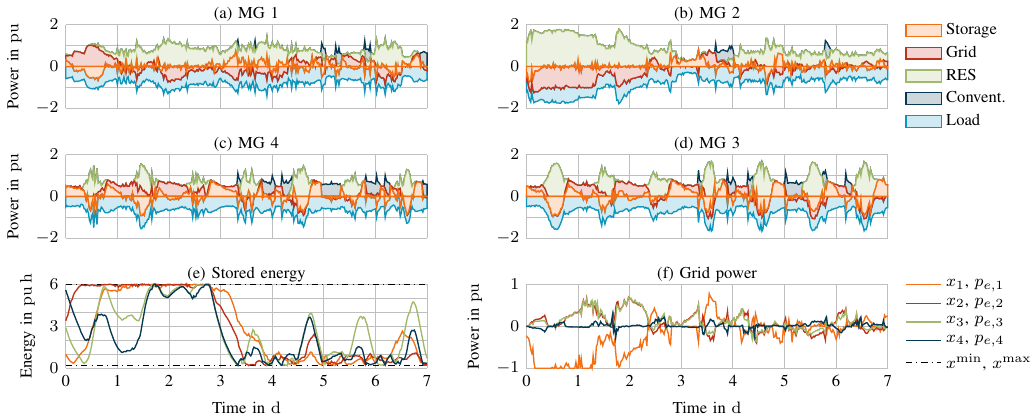}%
  \caption{Closed-loop simulation results obtained with \cref{algo:distributed}.}
  \label{fig:caseStudyResults}
\end{figure*}

Closed-loop simulations with three distributed approaches were performed:
\begin{enumerate*}[(i)]
\item \textit{Certainty equivalence \ac{mpc}}: This state-of-the-art approach is based on \cite[Algorithm~1]{HBR+2018}.
	The uncertain input is assumed to be given by the nominal forecast.
	Note that for comparability, switching costs were added to the approach from \cite{HBR+2018}.
\item \textit{Stochastic \ac{mpc}}: This approach is based on \cref{algo:distributed}.
	It employs forecasts in the form of scenario trees.
\item \textit{Prescient \ac{mpc}}: This approach is used as a reference and based on \cite[Algorithm~1]{HBR+2018}.
	Here, a hypothetical perfect forecast, given by the actual measured future load and available renewable infeed, is employed.
\end{enumerate*}
The setpoints found through them are passed to a model of the network of interconnected \acp{mg} which simulates the system behavior.
The case study was implemented in MATLAB~R2020b using
YALMIP (R20201001) \cite{Lof2004} and
Gurobi 9.1.0 as numerical solver.

\subsection{Closed-loop simulations}
\label{sec:closedLoop}

\Cref{fig:caseStudyResults} shows power and energy when using \cref{algo:distributed}.
The photovoltaic generators of \acp{mg}~3 and~4 come with a daily seasonality.
This is also reflected in power and energy of the storage units in these \acp{mg}.
Moreover, the power line between \acp{mg}~3 and~4 is hardly utilized, because of similar renewable infeed and load patterns.
The power provided by the wind turbines in \acp{mg}~1 and~2 does not come with notable seasonality.
\ac{mg}~2 exhibits a high renewable share on the first two days which is partly traded with other \acp{mg}.
From day 3 on, the available renewable infeed decreases.
Up to this time, all demand could be met without conventional generation.
After day 3, conventional units need to switch on at times to meet the demand.
However, their infeed remains small compared to \acp{res}.

\begin{figure}[b]
	\centering
	\includetikz{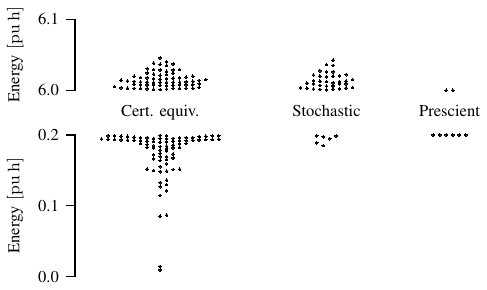}%
	\caption{Energy values outside of desired range. Motivated by the precision of the numerical solver, values closer than $10^{-5}$ to the interval were omitted.}
	\label{fig:violationsEnergy}
\end{figure}

\subsubsection{Violation of operating bounds}
\label{subsec:violations}

The storage units' soft bounds ${{x}}_{i}^{\min} = \unit[0.2]{pu\,h}$ and ${x}_{i}^{\max} = \unit[6]{pu\,h}$ for $i\in \N_{[1,4]}$ mark a desired range of operation.
\Cref{fig:violationsEnergy} shows a comparison with respect to this aspect.
We can see that the stochastic \ac{mpc} causes a smaller number of less extreme values outside of the desired interval compared to certainty equivalence \ac{mpc}.
Moreover, violations of power limits could only be observed with the certainty equivalence \acs{mpc}.
In total, 52 power limit violations (mean value: $\unit[0.012]{pu}$, maximum: $\unit[0.063]{pu}$) occurred.
Unlike the desired range of energy, these violations seriously jeopardize a safe operation.
The stochastic \acs{mpc} did not cause any power limit violation, which hints at improved safety when using \cref{algo:distributed}.

\subsubsection{Closed-loop performance}
\label{sec:costs}

\begin{table}[b]
	\centering
	\caption{Accumulated closed-loop simulation results}
	\label{tab:accumulatedResults}
	\begin{tabular}{lrrrr}
	\toprule
	 & & \parbox{11mm}{\centering Certaint\-y equival.} & \parbox{12mm}{\centering Stochasti\-c (Alg. \ref{algo:distributed})} & Prescient \\
	\midrule
	\multicolumn{2}{l}{Renewable energy in \unit{pu\,h}} & 318.1 & 333.8 & 334.0\\
	\multicolumn{2}{l}{Conventional energy in \unit{pu\,h}} & 61.0 & 45.2 & 45.9\\
	\multicolumn{2}{l}{No. of switching actions} & 46\phantom{.0} & 40\phantom{.0} & 29\phantom{.0} \\
	\midrule
	Costs \hspace{1em} & \acs{mg} 1 & 1\,652.8 & 1\,005.9 & 786.8\\
				& \acs{mg} 2 & 1\,366.1 & 774.2 & 607.7\\
				& \acs{mg} 3 & 2\,000.2 & 1086.1 & 1\,003.8\\
				& \acs{mg} 4 & 1\,755.3 & 1154.3 & 1\,109.9\\
				& Transmission & 21.7 & 19.0 & 18.7\\
	\cmidrule{2-5}
				& Sum & 6\,796.1 & 4\,039.5 & 3\,527.0\\
	\bottomrule
	\end{tabular}
\end{table}

\Cref{tab:accumulatedResults} contains the total generation of renewable and conventional units over the simulation horizon.
Moreover, accumulated closed-loop operating costs of all \acp{mg}, transmission costs and overall costs are summarized.
As expected, the prescient controller achieves the best results.
The certainty equivalence \ac{mpc} yields the worst costs and the stochastic \ac{mpc} takes the middle position.
Compared to the certainty equivalence case, renewable infeed of the stochastic approach is much higher and comes very close to the prescient case.
Moreover, the stochastic \ac{mpc} could reduce the share of conventional energy compared to the certainty equivalence \ac{mpc} by \unit[25.9]{\%}, the number of switching actions by \unit[13]{\%} and the overall cost by \unit[40]{\%}.

\subsection{Computational properties}
\label{sec:numeric}

\subsubsection{Solve times}
\label{sec:solvetime}

In what follows, we will discuss the accumulated solve times for each simulation step, i.e., each execution of \cref{algo:distributed} and \cite[Algorithm~1]{HBR+2018}.
The simulations were executed on a computer with an
Intel\textsuperscript \textregistered Xeon\textsuperscript \textregistered E5-1620 v2 processor @3.70 GHz with 32 GB RAM.
The certainty equivalence and the prescient \acs{mpc} yield very similar solve times (mean: $\unit[1]{s}$, maximum: $\unit[11]{s}$).
The stochastic \ac{mpc} requires a multiple of this (mean: $\unit[9]{s}$, maximum: $\unit[223]{s}$), which is noncritical considering a sampling interval of $\unit[30]{min.}$.
Compared to \cref{prob:mip} (mean: $\unit[127]{s}$), the solve times of \cref{algo:distributed} are significantly reduced and only come at the cost of a small suboptimality gap (mean: $\unit[0.3]{\%}$).

\subsubsection{ADMM convergence}

\begin{figure}
\centering
	\vspace{0.5em}
\includetikz{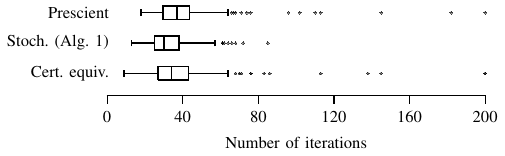}
\caption{Boxplots of number of \ac{admm} iterations for different approaches.}
\label{fig:admmIterations}
\end{figure}

\cref{fig:admmIterations} shows boxplots that illustrate the number of iterations required by the different algorithms when considering a tolerance of $\epsilon = 10^{-4}$ and $l_{\max} = 200$.
We can see that the median of all cases is below $40$.
\Cref{algo:distributed} comes with the lowest median which is as little as $30$ iterations and a maximum of $85$ iterations.
A deeper look into the simulation reveals that
outliers with a very high number of iterations seem to occur when all four \acp{mg} have similar conditions, i.e.,
either a surplus of \acp{res} and full storage units or
empty storage units and little available renewable infeed.
However, for a sampling time $\unit[30]{min.}$, $200$ iterations at maximum appear tolerable.

%% file: conclusions.tex

\section{Conclusions}
\label{sec:conclusions}

In this work, a scenario-based stochastic \acs{mpc} scheme for the operation of interconnected \acp{mg} was presented.
Based on a central \ac{mpc} formulation, a distributed algorithm that employs the \ac{admm} was developed.
The algorithm reflects the hierarchical structure in the network of interconnected \ac{mg}:
local controllers are in charge of individual \acp{mg}, while a central entity is in charge of the transmission grid.
In closed-loop simulations, the novel approach outperformed the certainty equivalence one concerning the number of constraint violations and the costs.
Moreover, the algorithm converges sufficiently fast for operation control.

Future work concerns simulations with a larger number of \acp{mg} and theoretical analyses concerning the scalability of the approach.
Moreover, suboptimality, persistent feasibility and privacy shall be further investigated.